*Perspective*

# Some Multifaceted Aspects of Mathematical Physics, Our Common Denominator with Elliott Lieb [†]

Daniel Sternheimer [1,2,‡]

1 Department of Mathematics, Rikkyo University, Tokyo171-8501, Japan; dasternh@gmail.com
2 Institut de Mathématiques de Bourgogne, 21078 Dijon, France
† Dedicated to our friend Elliott Lieb on the occasion of the ninetieth anniversary of his birth.
‡ Honorary Professor, St Petersburg University (Russia) and Member of the Board of Governors, Ben Gurion University (Israel).

**Abstract:** Mathematical physics has many facets, of which we shall briefly give a (very partial) description, centered around those of main interest for Elliott and us (Moshe Flato and I), and around the seminal scientific and personal interactions that developed between us since the sixties until Moshe's untimely death in 1998. These aspects still influence my scientific activity and my life. They also had as a corollary a variety of "parascientific activities", in particular, the foundation of IAMP (the International Association of Mathematical Physics) and of the journal LMP (Letters in Mathematical Physics), both of which were strongly impacted by Elliott, and Elliott's long insistence that publishers do not demand "copyright transfer" as a precondition for publication but are satisfied with a "consent to publish", which is increasingly becoming standard. This article being mainly a testimony to the huge scientific impacts of Elliott and also of Moshe, their intertwined aspects constitute the core of the present contribution. The last part deals briefly with metaphysical and metamathematical considerations related to axioms.

**Keywords:** mathematical physics; Elliott Lieb; the International Association of Mathematical Physics (IAMP, history and development); Letters in Mathematical Physics (LMP); history; God as an (optional) axiom; metamathematics; symmetries; particle physics

**MSC:** 01A99; 01A85; 01A80; 00A79; 00A30; 03F40

## 1. Some History and Related Material

### 1.1. The Context of Our First Interactions with Elliott Lieb

**a.** *Moshe Flato and I* first met at HUJI (the Hebrew University of Jerusalem) in 1958–1961, when we were M.Sc. students, he of Giulio (Joel) Racah and I of Shmuel Agmon, whose original lecture on the theory of distributions we both attended in 1958–1959. Agmon turned 100 in February 2022 and still lives in Jerusalem, where I visited him in early June; on that occasion we discovered that his French D.Sc. thesis had been edited for French language by a common friend (the wife of a leading French astrophysicist), whose father was a friend of my mother and had found for me a room in the apartment of a neighbor when I was student at HUJI. That year, I also visited Moshe's home in Tel Aviv to borrow the excellent notes he had taken at the lecture on hydrodynamics, given the year before by Abraham Robinson, who was then back at HUJI and is better known at present as the founder of "nonstandard analysis". (After getting a B.Sc. at the Hebrew University and, at the beginning of WWII, being stuck in France with a German passport, Robinson managed to reach London and became an expert on the airfoils used in the wings of fighter planes (essentially, conformal mapping) which he taught to himself while serving with the Free French Air Force). The lecture on hydrodynamics was a prerequisite for the "graduate seminar" on aerodynamics where Robinson had given me a talk and which I was attending at HUJI after immigrating to a kibbutz in Israel in October 1958 with a "Licence" (B.Sc.) from Lyon in France.

Quite naturally I was asked (by a relative of Moshe in Lyon who was my mother's best friend and her deputy in the local branch of the Women's International Zionist Organization, it is a small world) to help Moshe Flato after he arrived in Paris in October 1963, before defending the Ph.D. on group theory in nuclear physics, which he had prepared under Racah. Our interaction became much stronger, and its impact is felt after his untimely death in 1998. In Paris, since October 1961, I was "attaché de recherche" at CNRS, with Szolem Mandelbrojt (who had been Agmon's D.Sc. adviser) as adviser and his former student Jean-Pierre Kahane as deputy.

Incidentally, I had met Kahane at the first scientific meeting I attended, an "International Symposium on linear spaces", 5–12 July 1960 in Jerusalem. Much later, Lech Maligranda (a former student of Orlicz, who appears in the photo) published on arXiv [1] a very interesting historical note consisting of the photo below, accompanied by a brief description of the meeting. It shows the participants, most of them famous mathematicians (including Agmon); I still remember many of them. For the benefit of the reader, I insert in this paper, a semi-final version of the photo that was sent to me by Maligranda, with a request to help identify some of the few names that were not handwritten by him on the Figure 1.

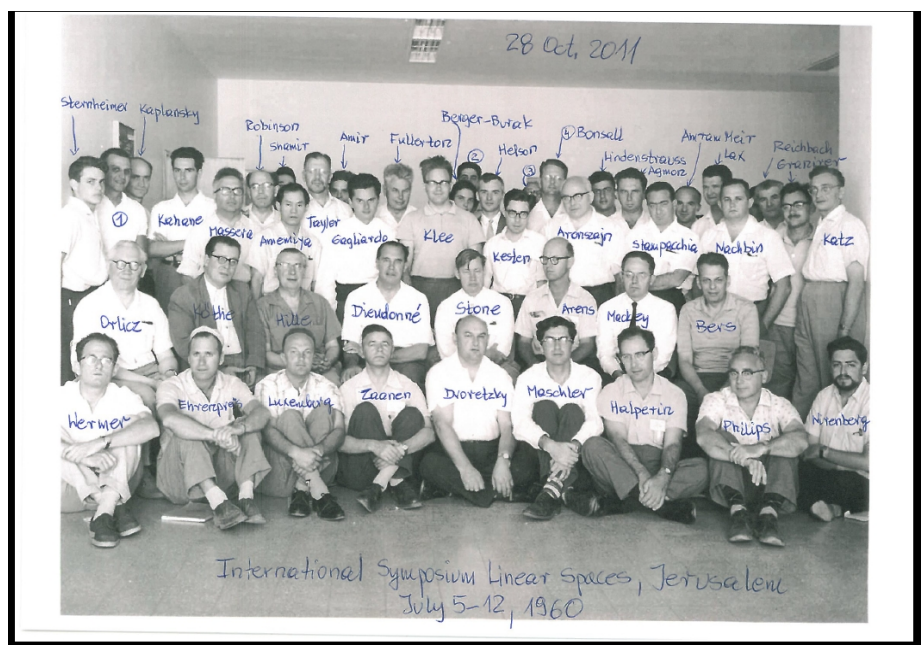

**Figure 1.** Jerusalem 1960, Linear Spaces Symposium.

In 1964, Moshe and I started to work on applications of group theory to physics (a major topic for Wigner and Racah), in particular, to particle physics. That evolved into "team work" (see, e.g., [2]), which lasted 34 years, until his death, and which I have been developing ever since.

**b.** The first *interactions* we (Moshe Flato and I) had with Elliott, beyond those which happen normally between scientists, date (if I remember well) to the mid seventies, after Elliott arrived at Princeton. It seems to me important to describe briefly how (and in which context) this happened. However, before this, I want to extend my most heartfelt congratulations to Elliott for being awarded (in January 2022) the highest honor bestowed by the American Physical Society, the APS Medal for Exceptional Achievement in Research for *"major contributions to theoretical physics through obtaining exact solutions to important physical problems, which have impacted condensed matter physics, quantum information, statistical mechanics and atomic physics"*. (The Medal was awarded for the first time in 2016 to Edward Witten.) And for being awarded the 2022 Gauss Prize at the International Congress of Mathematicians *"for deep mathematical contributions of exceptional breadth which have shaped* *the fields of quantum mechanics, statistical mechanics, computational chemistry, and quantum information theory.*" (Since 2006 the Gauss Prize is awarded every four years at the IMU Congress to honor scientists whose mathematical research has had an impact outside mathematics.)

Flato and I have been visiting Princeton, frequently for the time, since our first visit from France to the US in 1966, at the invitation of Eugene Wigner. Establishing a connection with Wigner was a natural step for Moshe. At Princeton, Wigner held the Higgins Professorship in physics, a named chair, which (after Wigner's retirement) Elliott inherited. The chair, in Wigner's time, was not subjected to regular teaching duties. However, when Elliott received it, and was asked by his colleagues about retirement, he declared very honestly that he does not intend to retire. Since he had only turned 43 in July 1975, the "perk" in teaching duties (not the named chair) became limited to 5 years! (Elliott transferred to Emeritus status only 42 years later.)

Already in the seventies, Moshe and I had experience in a wide variety of areas of mathematics (both applied and pure) and physics. This certainly helped us to then become friends with Elliott after his arrival at Princeton. A main factor for that "chemistry" was also that, though specializing in different aspects of mathematical physics, we shared a love for Science in all its aspects, especially in mathematics, physics and their interrelation.

Moshe's uncompromising attitude made a strong impact in the French scientific community immediately after his arrival in October 1963 (as a theoretical physicist and for some reasons in the group of Louis de Broglie). Eventually (and quite naturally), he found "scientific asylum" in the (more open) French mathematical school.

As mentioned above, I was, since 1961, a member of CNRS and a D.Sc. student in complex analysis, participating in the vivid mathematical life at Institut Henri Poincaré. That included, in 1963–1964, participating (together with many who eventually became leading French mathematicians in our generation) in the Cartan–Schwartz seminar [3] on the (then new) important Atiyah–Singer index theorem, which was held in parallel with another seminar on the same topic at Princeton organized by Richard Palais [4], exchanging by airmail roneotyped documents. My share (Exp. No. 22, 23p.) in the "Schwartz" side of it consisted of two talks on the multiplicative property of the analytic index, a crucial and nontrivial step needed for "dimensional reduction" (to varieties of dimension 2 and 1, and to the Laplacian and Dirac operators, respectively). My interests shifted to group theory in (particle) physics in 1965 because of my increasingly close collaboration with Flato.

*1.2. Our Mathematical Physics around the Seventies*

The main topics treated by Moshe and I in the 1970s dealt with two main intertwined aspects of physics, in general, and its mathematical formulation, in particular. They are, on the one hand, the importance (and an original use) of symmetries, especially in particle physics and in connection with relativity (a natural development of our work in the 1960s). That included the use of the AdS (Anti de Sitter) deformation of the Poincaré (inhomogeneous Lorentz) group and of its two special (most degenerate) representations, discovered in 1963 by Wigner's brother-in-law Dirac [5], who called them "singletons" and which we called Di and Rac (on the basis of Dirac's Bra and Ket); we used them, in particular, to interpret (also dynamically [6]) the photon as a composite of two singletons.

On the other hand, an original mathematical interpretation of quantization as a deformation of classical theories, now widely called "deformation quantization" [7,8]. The latter aspect relates to (and if I may say, constitutes a conceptual basis for quantum physics, which has been, in a variety of ways, a leitmotif in Elliott's huge scientific production.

The "chemistry" that immediately developed with Elliott was increased by the fact that some of Elliott's first papers have a common background with Moshe's first interests in physics, in particular with Moshe's 1960 M.Sc. Thesis with Racah on ionic energy levels in crystals [9,10].

1.2.1. A Couple of Short Explanations

*a.* As was observed already in 1964–1965, the Poincaré symmetry group of special relativity $SO(3,1) \cdot \mathbb{R}^4$ can be viewed as a deformation (in the sense that had been defined then by Gerstenhaber [11]) of the Euclidean symmetry group $SO(4) \cdot \mathbb{R}^4$ of Newtonian mechanics. The mathematically precise notion of deformation of groups and algebras is, in a way, an inverse operation to the "physical" notion of "contractions", which had been introduced in a more limited context in the 1950s, "in physics" by E.P. Wigner and E. Inonu [12], and by I.E. Segal in a side remark at the end of an article [13]. The latter notion has been studied and generalized by a number of people (for an informative more recent paper, see, e.g., [14]).

In those days, a natural question was asked, whether there is any connection between the experimentally guessed (by analogy with spectroscopic symmetries, of which Racah and Wigner had made ample use) unitary symmetries of elementary particles, especially the $SU(3)$ "internal" symmetry of the "eightfold way" (of Gell'mann and Neeman), and the Poincaré "external" symmetry. It would have made life easier for many at the time that there be no connection, but we objected [15] on mathematical grounds (dear to Elliott in other areas of physics) to "proofs" that the only connection possible was a direct product [16]. The "theorem" of Lochlainn O'Raifeartaigh was formulated at the Lie algebra level, where the proof is not correct because it implicitly assumes that there is a common dense domain of analytic vectors for all the generators of an algebra containing both symmetries. (Some later attempts by physicists to obtain similar results in a variety of contexts also contained implicit assumptions.) In fact, as it was formulated, the result is wrong, which we exemplified later with counterexamples. That did not prevent O'Raifeartaigh from becoming a friend. A direct product result was proved shortly afterward by Res Jost and, independently, by Irving Segal [17,18] but only in the more limited context of unitary representations of Lie groups.

*b.* In recent years, it dawned on me, based on the fact that *deformations of algebraic structures play a major role in physics*, that the above could be a false question. In a nutshell, we know that the Euclidean symmetry can be mathematically deformed to the Poincaré group of relativity by introducing a nonzero parameter $1/c$ (where $c$ is the velocity of light in vacuum), and that, in turn, the latter can be deformed into AdS (the anti De Sitter symmetry $SO(3,2)$) by introducing a tiny negative curvature in space-time. This permitted us to show that the photon may be considered as dynamically composite (of two Dirac singletons) and that so can the leptons [6,19].

Moreover, we know that Lie groups and algebras may be deformed into the so-called "quantum groups", but as Hopf algebras. The axioms for the latter were written in the 1940s, well before truly representative examples emerged from physics in the 1980s (in particular in Faddeev's Leningrad group in relation with quantum integrable systems). These "quantum groups" have an additional Hopf algebra structure, which makes them, in a way, analogous to Lie groups. They have been extensively studied in the past 40 years or so, and applied in various areas of physics. We also know since the 1990s that the latter "at root of unity" (i.e., when the deformation parameter is a root of unity) have finite-dimensional unitary representations, a property that was at the base of the introduction of compact "internal symmetries" to organize the experimentally discovered multiplets of elementary particles (not to mention the so-called "quarks", proposed already in 1964, but that is another story).

It is, thus, tempting to try and consider some form of quantum AdS as a candidate for these mysterious internal symmetries, even more so since it arises from relativity by deformations. That is what I have suggested in recent papers and talks (see, e.g., [20]).

1.2.2. The Context around This Contribution

As mentioned above, some of Elliott's works deal with topics that have a nonzero intersection with Moshe's early works (with Racah). More importantly, in most of his works and in ours, paying attention as much as possible to mathematical rigor is essential (which for Elliott may include finding the best constant in inequalities). Admittedly, a number of physicists do not understand the point in working so hard to prove mathematically results that have been "known" for a very long time using handwaving arguments, and have often been "confirmed experimentally". However, Elliott and us are convinced that achieving as much mathematical rigor as possible is of utmost importance, can prevent drawing erroneous conclusions (as shown also by our above-mentioned counterexamples) and may even lead to the discovery of new phenomena.

Furthermore, Elliott was never narrow minded and is (like us) interested in a large variety of scientific domains. His position at Princeton gave him superb occasions to satisfy his scientific curiosity. (For us, extensive traveling and inviting a wide spectrum of mathematicians and physicists, which Elliott also practiced a lot, achieved similar results). I remember that, during one of my visits to Princeton this century, he was amused by the fact that "stringies", as he called the many people working in and around the so-called string theory (which is more a framework as David Gross remarked), had become very excited by our works on singletons (see, e.g., [6]). In October 2015, he (together with Michael Aizenman) invited me to deliver a talk at the Princeton Mathematical Physics seminar (which meets at irregular intervals on Tuesdays). That talk was very well attended by leading people from the University and the Institute. My title (based in part on [20]) was (in obvious allusion to a popular paper by Wigner): "The reasonable effectiveness of mathematical deformation theory in physics, especially quantum mechanics and maybe elementary particle symmetries".

Our interactions with Elliott were not limited to the professional side, but here is not the place to expand on that. There were many occasions for interactions since he has always been an avid traveler, also for non scientific reasons (which did happen rarely to us and, the pandemic and age getting in the way, regretfully happens less to me in recent years).

## 2. Elliott and Us, the Science and Society Aspect

In the following, I shall therefore, concentrate on three main aspects of our interactions, which relate to physics and (the scientific) society: the birth of IAMP (the International Association of Mathematical Physics) in the 1970s; Elliott's impact on LMP (Letters in Mathematical Physics), the scientific journal initiated by Moshe, both in the mid 1970s; and his largely victorious battle with publishers on the Copyright issue.

LMP started with the relatively small D. Reidel publishing company, later included in the mathematical section of Kluwer. Eventually, LMP became affiliated at Kluwer with physics (for a variety of "corporate" reasons), after Elliott's time as one of the Editors. It remained there when the scientific part of Kluwer was merged by new owners with Springer into a new company, still named Springer, which is now part of the huge Springer–Nature. That is a typical example of the acquisition-merger trend, which pervaded also the scientific publications world. It has some benefits (of scale in particular) but it is potentially dangerous in our increasingly digitalized era. For example, what will happen to the few platforms that host most scientific publications, paid for by the work of dedicated scientists and the (large) contributions of their institutions, if the corporations who own and manage these platforms become bankrupt and the platforms become suddenly dark? It would take time for the scientific societies to get around the technical and legal problems involved, and in the meantime, most of the past scientific work will be available only in the libraries of a few institutions, bringing (especially in theoretical domains) research tools back at least half a century.

### 2.1. The Birth of IAMP

*(a). The prehistory of $M \cap \Phi$.* In April 1966, the CNRS organized in Gif-sur-Yvette an international conference titled "Extension du groupe de Poincaré aux symétries internes des particules élémentaires". Flato was an initiator and a co-organizer, together with Louis Michel and Jean-Pierre Vigier. At the conference dinner, he sat next to Gunnar Källén, an auxiliary member of the Nobel Committee for physics, who, like Flato, had (mildly speaking) a sharp tongue. (I was seated on the other side of Källén and remember well some of their exchanges, which are not for publication, even now when those involved are no longer with us). Between the two developed an immediate and strong empathy. A year later, our friend Gilbert Karpman was appointed scientific attaché of France in Stockholm and a Franco-Swedish conference on mathematical physics was planned with Flato and Källén. After Källén's death in the crash of the plane, he was piloting on his way to CERN in October 1968, that series of meetings (the first was held in Stockholm in December 1968) was given Källén's name. A second one was held in Paris in June 1970 and a third in Göteborg in June 1971. In December 1972, that became a franco-polish-swedish conference on fundamental problems in elementary particles physics, held in Warsaw just before the "International Conference on Mathematical Problems of Quantum Field Theory and Quantum Statistics" magnificently organized in Moscow by our friend Nikolay Nikolayevich (N.N.) Bogoliubov (with a ceremony at the Kremlin). The latter was eventually considered as the First IAMP Congress. That is where the symbol $M \cap \Phi$ was introduced.

It was at the International Congress of Mathematicians held in Moscow in August 1966, which Moshe and I attended as part of the (relatively large) French delegation, that we first met N.N. Bogolyubov and many other Soviet mathematicians with whom Moshe established friendly relations from the beginning. Among them were the young Ludwig Dmitriyevich Faddeev and the older Israïl Moyseyovich Gelfand, who invited Moshe to deliver a talk at his celebrated seminar. There was immediate empathy with them, facilitated by the fact that Moshe could speak quite fluently Russian with an almost native accent which his perfect ear had caught from his family (he had, however, to rely on me to read Cyrillic). At ICM66 N.N. (Nikolay Nikolayevich) invited Moshe and me to come after ICM66 to the Laboratory of Theoretical Physics (which he created and now bears his name) in the relatively new J.I.N.R. (Joint Institute for Nuclear Research) established in 1956 in the new "town of science" Dubna. Which we did, and we visited Dubna several times after that. I continued the tradition this century until the pandemic (and more) got in the way.

In December 1972, after the extension to Poland of our Gunnar Källén meetings, almost all participants of it continued to Moscow to $M \cap \Phi$. A cocasse anecdote occurred in Warsaw while we were there for Christmas 1972 at the home of our friend, the late Ryszard Rączka. He had direct phone connection to Moscow, so Moshe called I.M. Gelfand's home and the call started as follows: "Merry Christmas Israïl Moyseyovich" said Moshe, to which the latter replied "Merry Christmas Moysey Salomonovich".

*(b).* In March 1974, a continuation and extension of the 1972 franco-polish-swedish meeting was held in Warsaw, an International Symposium on Mathematical Physics (eventually considered as the second IAMP Congress). That is where we launched the ideas of both a mathematical physics society and a new scientific journal, of shorter publications, a somewhat mathematical physics analog of the Physical Reviews Letters.

*2.2. The Development of the Concept of IAMP*

Though a European, I had become an individual member of the APS (American Physical Society) in 1967, after our first visit to US in 1966 during which we visited at BNL (Brookhaven National Laboratory) my direct cousin Rudolph Sternheimer, who sponsored me. APS is, by nature, much larger than the AMS (American Mathematical Society) and organized differently. A European Physical Society (EPS) was created only in 1968, mainly as a federation of national societies (individual membership is possible, but relatively rare and powerless). The European Mathematical Society (EMS) was founded much later, in 1990. It is also a federation of about 60 national societies, but has many individual members electing representatives who participate in the general assemblies alongside with delegates of national societies.

After the creation of the EPS, it seemed natural to me, since mathematical physics was more developed in Europe, that a European Mathematical Physics society be created, but (like the APS) with mainly individual members, who could be coming from all over the world. That is why, together with Flato, we suggested just that in 1974 in Warsaw. There was an immediate reaction, in particular from our friends Elliott and Huzihiro Araki (from Japan, who also turned 90 in July 2022) who said: the idea is good but should by no means be restricted to Europe in its denomination. We obviously agreed. The development took some time, in particular because a number of leading mathematical physicists did not see the necessity and were afraid that Moshe would use it for personal enlargement, which had never been his intention.

Rudolf Haag, who, at first, had been reluctant to the idea, eventually became convinced that Moshe was not looking to use IAMP. That happened in 1975, at a meeting in Lausanne (managed by Marcel Guenin, who had played a role in the early EPS). He described the events in a historical paper [21] from which I extract the following part:

> *"An important development for mathematical physics taking shape in 1975 was the foundation of the International Association of Mathematical Physics. The creation of such an organization had been proposed for several years by Moshe Flato and pushed very vigorously by him against some opposing faction of scientists which included me. The controversy was in part due to lack of clarity about the objectives of the proposed organization but in part also due to personal animosities. Some of us had begun our scientific life before the great inflation in numbers at a time when the theoretical physics community was a rather tightly knit group, inspired by great masters like Lorentz, Planck, Einstein, Bohr, Sommerfeld, ... We did not see any need for a new organization outside the existing mathematical and physical societies and feared the spectre of a public relations oriented lobby engaged in fund raising for some pet projects. The somewhat flamboyant and aggressive manners of Moshe Flato had earned him quite a number of enemies. Res Jost had published* [17] *an unusually sharp reply to a criticism by Flato; Louis Michel had had some clashes with him; Daniel Kastler and myself were embarrassed and annoyed when at a party in Moscow, while we were talking with Bogolubov, Flato came up raising his glass, slapping Daniel on the shoulder exclaiming: "Don't you think Daniel, that we should see to it that Bogolubov gets the next Nobel Prize? But Flato had also friends who appreciated his unconventional ways and his generosity. A 1974 attempt to create the organization by an overwhelming vote of the participants at an international congress on mathematical physics in Warsaw failed, mainly because the Russian delegation was uncertain whether this was politically correct. So in Fall 1975 it was decided that a few representatives of the opposing groups should get together and settle the issue. We met in Lausanne. On one side there was Flato and Piron, on the other side Hepp and myself and, if I remember correctly, Borchers as a neutral witness. In the course of the discussion Flato succeeded in convincing me that he was not a bad guy and we ultimately agreed that the organization should be created, that the first president should be Walter Thirring and that in the executive board there should be no person who had played any role in the previous controversy. Thirring accepted the task and appointed a committee of four persons, consisting of Araki, Piron, Ruelle and Streater, to work out the statutes of the organization. Araki in his usual careful, conscientious way wrote the final version of the statutes, which were approved by the vote of the inscribed members in July 1977. Thus the organization could start its life".*

The above-mentioned paper by Haag was initiated by (and with) Araki and me, for a different purpose that was eventually discontinued by the publisher and the paper was published in the H (history) section of the EPS Journal. It basically reflects the events of that time, but, in my opinion, a number of statements could require "footnotes", especially concerning the description of events involving Flato and me. Most of that is another matter (see, e.g., [2]), but some precisions can be made, in addition to the fact that I came to Moshe with the idea of an Association, which he liked and adopted as his: we often operated like that, as a team. (The issue of the Nobel prize for Bogolubov, which did not happen for reasons related to the way the Soviet Union operated, is still restricted.)

As a "sabra", Moshe seemed quite extroverted with a flamboyant and often aggressive style. However, that was his way to hide shyness and humility, as people who got to know him well eventually realized, including Yvette Chassagne, the first woman to become "préfet" in France, after which she was chairperson of Union des Assurances de Paris (UAP) between 1983 and 1987, then the largest French insurance company. There, she asked Flato to establish a Scientific Council which included many VIPs, and a scientific prize, that was awarded by that Council, one of the first being given to Lieb in 1985.

Moshe did not leave anybody indifferent. Most of those who got to know him became friends, which includes Elliott of course. Moshe's extraordinary personality and ability to develop contacts in many areas and countries gave him a special status. Mutatis mutandi, similar things can be said of Elliott.

*2.3. The Prolonged Impact of Elliott on IAMP*

2.3.1. The Evolution of IAMP

As attested by Haag, Elliott's friend Walter Thirring played a crucial role in the effective start of IAMP and was the president during its first 3 years of existence (1976–1978), followed by Araki and Elliott, who was always in the background and is the first scientist to have been president twice (1982–1984 and 1997–1999). As a matter of fact, the indirect impact of Elliott on IAMP never ceased to be felt and remains a kind of watermark, though he himself would have preferred, and been more comfortable with, a greater diversity in the definition of "IAMP mathematical physics", which represents only a part of his wide scientfic interests. Elliott (private communications) is fully aware that this is not a healthy situation, and so are many members of the Executive Committees over the years. Serious efforts towards increasing diversity are made, but so far inner dynamics make these efforts insufficient. Both Elliott and I feel it would be good for the spectrum represented in the IAMP institutions and activities to be more inclusive.

2.3.2. A Perverse Effect of Democracy

As one can read in its statutes (which appear as a page in the IAMP site), IAMP is governed by an "Executive Committee" elected by a ballot of the General Assembly (in practice, electronically before an ICMP) for a term of three years (renewable only once) by its ordinary members who have paid their (modest) dues. The problem, which happens to various degrees in many democratic institutions, is that many mathematical physicists (an admittedly imprecise notion) do not pay dues or do not bother to vote for a variety of reasons, and that those who do vote represent only a fraction of the spectrum of mathematical physics, concentrated in parts of the spectrum. The net result is that the first twelve scientists appearing on the list after the vote increasingly reflect these parts, which contributes to discouraging people from other areas.

It is natural that scientists tend to vote for scientists they know, which often means people close to their fields. Not many have as wide a vision of Science as Elliott. Being at the origin of the creation of IAMP, many of its first members knew me well. That is how I was elected to its first two Executive Committees, and in that capacity, became also a member of the Commission on Mathematical Physics (C18) established by the IUPAP (the International Union of Pure and Applied Physics) in 1981. However, a perverse effect manifested itself with the passing of time. The Executive Committees tried to suggest that voters take into account the diversity of mathematical physics, which, in fact, might be better represented by the symbol $M \cup \Phi$, in order to make more clear that it includes also what we call "physical mathematics", i.e., mathematics motivated by physics. Incidentally, in the UK, for a long time, mathematical physics referred mainly to the study of partial differential equations.

The problem recently became more acute with the pandemic and the overwhelming use of talks via Zoom or similar, which made an increasing number of talks accessible to a large audience (if people have time for that). As stated on the IAMP site, a long list of talks can be found on the site researchseminars.org, where the official "One World" IAMP seminars are listed. That long list is not easy to use, and (at present) misses a number of important seminar series that give perspectives on wide areas. For instance, the Rutgers

Mathematical Physics seminars of Elliott's friend Joel Lebowitz who (at 91) manages to bring, week after week, review talks by leading speakers covering important sectors of mathematical physics and applied mathematics in a very broad sense.

*2.4. Elliott and LMP*

Concommitant to, but independant of, the creation of IAMP is the start of LMP (Letters in Mathematical Physics), which initially was meant to be a journal of important short contributions, being to CMP (Communications in Mathematical Physics) a kind of analog of what is PRL (the Physical Review Letters) is to the remainder of the Physical Review. It was also initiated by Moshe in Warsaw in 1974. The first contacts were made (with the help of Ryszard Rączka) with PWN (Polish Scientific Publishers), but it soon became clear that, given the delays that regular mail would bring (in particular due to political censorship in Poland), that was not a realistic option. However, through PWN, contacts were made with the (then small) D. Reidel Publishing Company, based in Dordrecht, whose owner Anton Reidel liked to come to Poland in order to buy icons.

That is how the first issue of LMP appeared at the end of 1975, with at first four Editors (Flato, Ryszard Rączka, Stan Ulam and Marcel Guenin) and a diversified Editorial Board. For a short while thereafter, in order to speed up publication, it was decided that proofs would not be sent to the authors unless they insisted on it. After a year or so, during which errors were introduced at the production level by an automated text editor (e.g., "conformal" becoming "conformational" in a title), a problem we all encounter until now with automatic corrections on many devices, that practice was abandoned. Eventually the publisher grew in size, becoming one of the two major Dutch publishing houses in Science, and hired for LMP a superb dedicated editor for the papers (a former British physics student named Richard Freeman). LMP also grew in size and so did its backlog, but the quality remained high. Contrarily to what some had feared, there was no competition with the "classic" journal Communications in Mathematical Physics (CMP, published by Springer in Heidelberg), rather complementarity. The structure of the journal evolved somewhat with time, but its way of operating remained.

Following our suggestion, Elliott joined LMP as one on the main Editors starting with volume 8 (issue 1), January 1984. That was announced in an Editorial signed by Moshe in the preceding issue. His first duty as Editor was to cosign with Moshe an obituary for Mark Kac, who had been a member of the Editorial Board of LMP. Elliott's impact was felt in many areas, both scientific and concerning the journal governance. He formally ceased being one of the Editors after volume 39 (issue 4), March 1997. When Flato died suddenly (at the age of 61) in November 1998, it fell on me to manage a smooth transition, which I did, with the help of the team he had built. Elliott's input has been very important in those circumstances, and he continues informally to be of help in many matters concerning the Journal.

*2.5. Copyright Transfer vs. Consent to Publish*

Before, during and after, our interactions around LMP, we became involved in the uphill battle Elliott has been waging, for years, with many publishers concerning the latters' demand of transfer of copyright for scientific texts that were published, mainly in journals. Without such a transfer, papers could not be published. (The case of books by one or few authors may be different). We supported Elliott's point and made it ours, including for our publications. With good reason, he considered that such a demand is outrageous, dealing with work performed by scientists (usually paid by academic institutions, mostly with public support and often outside of their obligations) and was the fruit of their minds, formed by years of studies. One of the arguments of the publishers, which we encountered both as authors and as editors, was that they are better equipped than scientists to defend the copyright (from possible plagiarism) than individual authors. In mathematical physics (and other areas), such an argument is largely hypothetical. We were sometimes told (orally) that there is a notion of "fair use" permitting authors to use (up to about 20% of) their own work in a later work. That would not be needed if only Consent to Publish is granted. However, even if Copyright is transfered, no sensible publisher would count words and sue an author if it determines that the use of previous material exceeds the interpretation of "fair use".

My conjecture is that a hidden reason for the demand of copyright transfer, across the board of the large spectrum of a publisher (without consideration of the fact that domains like mathematical physics can be special cases), is that some "smart" managers assumed that the amount of copyrights they have adds to the market value of the publisher, if and when (as has happened an increasing number of times) some financial institution wants to buy it. This argument was probably found fallacious, the consent to publish, which Elliott has been insisting on since many years, being sufficient. At first, a compromise was found, but not publicized for a long time, that authors may keep Copyright to themselves and give the publisher only "Consent to Publish" of the text in a given journal (a kind of "Lieb exception"). That is now increasingly the case of most of our publishers. Various juridical formulations are given to that, depending on the lawyers who play an increasing role with publishers. Here also, Elliott's vision and persistence have been rewarded.

### 3. Metaphysical—"Metamathematical" Remarks
*3.1. Preparatory Introduction*

On the web site of the journal "Axioms", its aims start with the sentence *"Axioms (ISSN 2075-1680) is an international, open access journal which provides an advanced forum for studies related to axioms"*. The following remarks are motivated by that sentence, though they are largely disconnected from the remainder of that contribution, and only vaguely related to the very general definition of the journal as "a journal of mathematics, mathematical logic and mathematical physics".

These remarks were, in part, triggered by a question that is often asked in various forms. That kind of question was recently (on 22 January 2020 to be precise) asked to me by the most tolerant rabbi of one of the two branches of Chabad in Tokyo: whether I believe in a God or at least in a kind of supreme being.

My answer started with the 1931 incompletemess theorems of Gödel in axiomatic set theory (more precisely, about natural numbers). The interested reader can find online precise formulations of these theorems. For the purpose of the present remarks, we can be satisfied with the somewhat vague statement that the (second) incompleteness theorem, an extension of the first, shows that such a system (of axioms) cannot demonstrate its own consistency.

My answer to the Chabad rabbi continued with a far reaching extension of the above, which, in a way, is the basis of my "Weltanschauung" (comprehensive conception of the world), according to which the notion of a God, or supreme being, is an axiom that may, or not, be added to a Weltanschauung. That attitude may be called "atheism" if one takes at its face value the prefix "a", in contradistinction with what Proudhon and others called "anti-theism" (a notion which appeared already in 1788 in the Oxford English dictionary).

For me the notion of a God (or supreme being) is an axiom which may, but need not, be added to whatever "axioms" you choose to govern your life. Both alternatives are self-consistent. It is then a matter of choice whether one chooses to add such an axiom to one's own system of axioms. That is, in a way, similar to the so-called "axiom of choice" in the most commonly accepted (Zermelo–Fraenkel) set theory, which is at the basis of most parts of modern mathematics. The analogy (imperfect as all analogies are) goes even further, because most humans use, in their everyday language, an implicit reference to some God. For example, I have often heard Japanese (most of whom have a very relaxed attitude toward religions) say "I *play* to God" (Asians often do not distinguish between r and l) merely to express some wish.

Incidentally, when in Fall 1958, I arrived from France at the Hebrew University (HUJI) as a (third year mathematics) student, Avraham (Halevi, born Adolph) Fraenkel (1891–1965) had just retired. My mother met him there in 1957, which was instrumental in her decision to let me immigrate to Israel in 1958 (at that time I was only 20 and still a minor by French law). Fraenkel was succeeded in his chair by Abraham Robinson (mentioned in the historical introduction above). There are many anecdotes related to Fraenkel from the time when he was teaching at HUJI (one of which I used in a small humoristic column that got published in "𝔏𝔢 𝔐𝔬𝔫𝔡𝔢" in the late 1950s), but that would take us too far off.

*3.2. God as an Axiom—And Some Corollaries*

My purpose here is not to add a tiny fraction, which some may consider as blasphematory, to the many (often millenary) discussions around the notions of God (sometimes written G*d, especially by "Haredi" Jews, a term derived from the Biblical verb *hared*, which appears in the Book of Isaiah and is translated as "[one who] trembles" at the word of God) and around the various religions derived therefrom.

My late mother, whom I respect immensely (fifth commandment of Judaism), used to tell me that a main difference between Judaism and all other religions is that every Jew has a direct line to God (if I may add, whatever that means), while in all other religions, there must be some intermediary (whatever name is given to that notion). The problem is that in most "sects" (including of Jews) the members feel the need to follow some guide. If that is their free choice, let it be.

My very personal (rather iconoclastic and deliberately carricatural) definition of religion is that it is a sociological phenomenon of sado-masochistic nature, in which those I will call here (for the purpose of cartoon)"sadists", avatars of those who claim to speak in the name of their god (and would like to impose their view on others) fix rules, which "masochists" enjoy imposing upon themselves. In "Le médecin malgré lui" by Molière (act 1, scene 2) a bigot (Robert) attempts to interfere in the life of a couple of servants (Martine and Sganarelle) because the latter beats Martine, only to be answered by her: "I want him to beat me".

Coming back to axioms, as stated above, my personal claim is that, by an audacious generalization of Gödel's incompleteness theorem, the notion of God is an axiom which may, but need not, be added to whatever set of axioms one uses as a guide. Many prefer to add such an axiom, and to follow one of the manifold variations on that theme that have been added over time. In many modern societies one is free not to do so, but that does not mean that one is "off the hook", on the contrary. When one lives in a society (which is the case of almost all humankind), there are limitations to the freedom of each individual person, when (and hopefully only when) it infringes upon the freedom of others. These limitations, or at least many of them, may be easier to accept when it is claimed that they come from some supreme being. However, claiming that the existence of such an axiom is not your problem does not free you from its corollaries.

On the contrary, it renders each and every one responsible to make life in society as harmonious as possible, without having to fear a "Père Fouettard", the character who supposedly accompanies Saint Nicholas on his rounds during Saint Nicholas Day (6 December). "Vaste programme" as is supposed to have said, in a very different context, Général De Gaulle after seeing on the first Jeep who entered Paris on 24 August 1944 an inscription "death to idiots" ("mort aux cons"). Indeed, neither alternative is easy to achieve, but these are not "mission impossible". Both Theodor Herzl and Walt Disney, and many others, have quotes going in similar directions, e.g., "If you will it, it is no dream" and "If you can dream it, you can do it. Remember it all started with a mouse".

**Funding:** This research received no external funding.

**Institutional Review Board Statement:** Not applicable.

**Data Availability Statement:** Not applicable.

**Conflicts of Interest:** The author declares no conflict of interest.